\documentclass[11pt]{amsart}

\usepackage{amsmath}
\usepackage{amssymb}
\usepackage{amscd}
\usepackage[all]{xy}

\makeatletter
\@addtoreset{equation}{section}
\makeatother

\newtheorem{theorem}{Theorem}[section]

\newtheorem{lemma}[theorem]{Lemma}
\newtheorem{corollary}[theorem]{Corollary}

\def\tg{\mathcal{T}(G)}
\def\tnqg{\mathcal{T}_{nq}(G)}

\def\mapg{G^{ap}}
\def\tauap{\tau_{ap}}
\def\etaap{\eta_{ap}}
\def\taunq{\tau_{nq}}
\def\Ree{\mathbb{R}}
\def\eps{\varepsilon}
\def\fQ{\mathcal{Q}}
\def\fS{\mathcal{S}}

\def\hdg{\mathfrak{HD}(G)}
\def\spineg{\mathrm{A}^*(G)}

\def\endpf{{\hfill$\square$\medskip}}
\def\proof{{\noindent{\bf Proof.}\thickspace}}
\def\comp{\raisebox{.2ex}{${\scriptstyle\circ}$}}

\begin{document}


\title[The spine]{The spine of a Fourier-Stieltjes algebra:  corrigenda}

\author{Monica Ilie and Nico Spronk}


\maketitle

\footnote{{\it Date}: \today.


We like to thank NSERC for partial funding of our work.}

It has come to the authors' attention that there are several 
errors in our paper \cite{ilies2}.  Fortunately,
these errors are correctable. 
Some of these errors carry on to modest, though mainly
cosmetic, errors in follow-up papers \cite{rundes} and \cite{ilies1}.  

We are grateful to El\c{c}im Elgun
for pointing out the gap in the proof of \cite[Theorem 2.2]{ilies2}
and the error in the proof of \cite[Theorem 5.1]{ilies2}.
We are also grateful to Pekka Salmi 
for pointing out the flaw in the statement of \cite[Theorem 4.2]{ilies2}
and suggesting its correct form.

We appeal to \cite{ilies2} for pertinent notation and terminology.

\section{On non-quotient locally precompact topologies}

The following should replace \cite[Theorem 2.2]{ilies2}.

\begin{theorem}\label{theo:tnqg}
Let $\tau\in\tg$.  Then the following hold.

{\bf (i)} There exists a unique $\taunq$ in $\tnqg$
such that $\tau$ is a quotient of $\taunq$.

{\bf (ii)} If $G$ is abelian, then $\taunq=\tau\vee\tauap$.
\end{theorem}

In \cite[Theorem 2.2]{ilies1} it is claimed that $\taunq=\tau\vee\tauap$
for general locally compact groups $G$ and $\tau$ in $\tg$.  This is
false.  For example let $G=\mathrm{SL}_2(\Ree)$.  We have
that, on $G$, $\tauap$ is the trivial topology $\eps=\{\varnothing,G\}$.
Now if $q:G\to G/\{-I,I\}$ is the quotient map
and $\tau=q^{-1}(\tau_{G/\{-I,I\}})$ is the coarsest topology making $q$ 
continuous, then
we have that  $\taunq=\tau_G$ whereas $\tau\vee\tauap=\tau\subsetneq\tau_G$.

The proof of part (ii) proceeds exactly as does the proof of 
\cite[Theorem 2.2]{ilies2}.  In particular, in the second paragraph of that 
proof, \cite[Lemma 2.3]{ilies2}
may be used to show that $s\mapsto\bigl(\eta^{\tau_1}_\tau(s),\etaap^{\tau_1})
:G_{\tau_1}\to G_\tau\times \mapg_{\tau_1}$ is a bicontinuous isomorphism,
since $G_{\tau_1}$, being abelian, is maximally almost periodic.

The proof of part (i), however, demands more care.  We fix
$\tau_0$ in $\tg$ and let
\begin{equation}\label{eq:quetau}
\fQ_{\tau_0}=\{\tau\in\tg : \tau_0\text{ is a quotient of }\tau\}.
\end{equation}

\begin{lemma}\label{ref:thelemma}
{\bf (i)} If $\tau_1,\tau_2\in\fQ_{\tau_0}$ then $\tau_1\vee\tau_2\in
\fQ_{\tau_0}$ as well.

{\bf (ii)} If $\tau_1,\tau_2\in\tg$ satisfy
$\tau_0\subseteq\tau_1\subseteq\tau_2$ and $\tau_2\in\fQ_{\tau_0}$,
then $\tau_1$ is a quotient of $\tau_2$.
\end{lemma}

\proof {\bf (i)}  For $j=1,2$ we let $K^j_0=\ker\eta^{\tau_j}_{\tau_0}$.
Our assumptions provide that for $j=1,2$, $K^j_0$ is compact and
$G_{\tau_j}/K^j_0\cong G_{\tau_0}$.  We identify $G_{\tau_1\vee\tau_2}$
as a subgroup of $G_{\tau_1}\times G_{\tau_2}$ and define
$K=G_{\tau_1\vee\tau_2}\cap(K^1_0\times K^2_0)$.  We let 
$q:G_{\tau_1\vee\tau_2}\to G_{\tau_1\vee\tau_2}/K$ denote the quotient map.
We obtain the following commutative diagram
\[
\xymatrix{
G_{\tau_1\vee\tau_2} \ar[d]_{q} \ar@{^{(}->}[r] 
& G_{\tau_1}\times G_{\tau_2} 
\ar[dr]^{\eta^{\tau_1}_{\tau_0}\times\eta^{\tau_2}_{\tau_0} } \\
G_{\tau_1\vee\tau_2}/K \ar@{^{(}->}[r]
& G_{\tau_1}\times G_{\tau_2}/(K^1_0\times K^2_0)
\cong G_{\tau_1}/K^1_0\times G_{\tau_2}/K^2_0
\ar@{=}[r]^<<<<<{\sim}
& G_{\tau_0}\times G_{\tau_0}.
}
\]
Identifying $G_{\tau_0}$ with the diagonal subgroup of
$G_{\tau_0}\times G_{\tau_0}$, this diagram
shows that the map $\eta^{\tau_1\vee\tau_2}_{\tau_0}
=\eta^{\tau_1}_{\tau_0}\times\eta^{\tau_2}_{\tau_0}|_{G_{\tau_1\vee\tau_2}}$
is a proper map.

{\bf (ii)}  We recall that our assumptions give the following commuting diagram
\begin{equation}\label{eq:firstdiagram}
\xymatrix{
G_{\tau_2} \ar[dr]_{\eta^{\tau_2}_{\tau_0}}
\ar[rr]^{\eta^{\tau_2}_{\tau_1}} 
& & G_{\tau_1} \ar[d]^{\eta^{\tau_1}_{\tau_0}} \\
& G_{\tau_2}/K^2_0 \ar@{=}[r]^<<<<<{\sim} &  G_{\tau_0} 
}
\end{equation}
where $K^2_0=\ker\eta^{\tau_2}_{\tau_0}$ is compact.  Then
$K^2_1=\ker\eta^{\tau_2}_{\tau_1}\subset K^2_0$, and is thus compact.
We let $q^2_1:G_{\tau_2}\to G_{\tau_2}/K^2_1$ be the quotient map and
$\tau=(q^2_1\comp\eta_{\tau_2})^{-1}(\tau_{G_{\tau_2}/K^2_1})$,
so that $G_\tau\cong G_{\tau_2}/K^1_1$, $\tau\supseteq\tau_1$
and $\eta^\tau_{\tau_1}$ is injective.
The commuting diagram (\ref{eq:firstdiagram}) and
the first isomorphism theorem give the commuting diagram
\begin{equation}\label{eq:seconddiagram}
\xymatrix{
G_{\tau_2} \ar[drr]_{\eta^{\tau_2}_{\tau_0}} \ar[r]^{\eta^{\tau_2}_{\tau}}
& G_\tau \ar[dr]^{q} \ar[rr]^{\eta^{\tau}_{\tau_1}} 
& & G_{\tau_1} \ar[d]^{\eta^{\tau_1}_{\tau_0}} \\
& & G_{\tau}/K \ar@{=}[r]^<<<<<{\sim} &  G_{\tau_0} 
}
\end{equation}
where $K=K^2_0/K^2_1$ and $q:G_\tau\to G_{\tau}/K$ is the quotient map. 

Since $\eta^\tau_{\tau_1}$ is injective, it suffices to prove
that it is open for $\eta^{\tau_2}_{\tau_1}$ to be a quotient map,
i.e.\ we obtain that $\tau_1=\tau$ where $\eta^{\tau_2}_{\tau}$
is a quotient map.  To this end, let $U\subset G_\tau$ be relatively
compact open set.  Then $UK$ is also relatively compact and open.
Hence $\eta^\tau_{\tau_1}(\overline{UK})$ is closed and equal to
$\overline{\eta^\tau_{\tau_1}(UK)}$, and furthermore
\begin{equation}\label{eq:homeomorphism}
\eta^\tau_{\tau_1}|_{\overline{UK}}:\overline{UK}\to
\overline{\eta^\tau_{\tau_1}(UK)}\text{ is a homeomorphism.}
\end{equation}
The commuting diagram (\ref{eq:seconddiagram}) tells us that
\[
\eta^\tau_{\tau_1}(UK)
=(\eta^{\tau_1}_{\tau_0})^{-1}\bigl(q(U)\bigr)
\]
which is open in $G_{\tau_1}$.  Hence by (\ref{eq:homeomorphism}),
\[
\eta^\tau_{\tau_1}|_{UK}:UK\to\eta^\tau_{\tau_1}(UK)
\]
is a homeomorphism onto an open subset so 
$\eta^\tau_{\tau_1}(U)$ is open.  \endpf

\noindent {\bf Proof of Theorem \ref{theo:tnqg} (i).}
We first note that $\fQ_\tau$, as defined in (\ref{eq:quetau}),
is a directed system:  if $\tau_1,\tau_2\in\fQ_\tau$ then
$\tau_1\vee\tau_2\in\fQ_\tau$ by Lemma \ref{ref:thelemma} (i).
Moreover, it follows Lemma \ref{ref:thelemma} (ii) that
if $\tau_1,\tau_2\in\fQ_\tau$ with $\tau_1\subseteq\tau_2$, then
$\eta^{\tau_2}_{\tau_1}$
is a proper map.  Hence the inverse mapping system
\[
\{G_{\tau'},\eta^{\tau_2}_{\tau_1}:\tau'\in \fQ_\tau,\tau_1\subset\tau_2
\text{ in }\fQ_\tau\}
\]
gives rise to the projective limit 
\begin{align*}
G_{\fQ_\tau}
=\underset{{\tau'\in\fQ_\tau}}{\underset{\longleftarrow}{\lim}}G_{\tau'}
&=\left\{(s_{\tau'})\in\prod_{\tau'\in\fQ_\tau}G_{\tau'}:
\eta^{\tau_2}_{\tau_1}(s_{\tau_2})=s_{\tau_1}\text{ if }\tau_1\subseteq
\tau_2\text{ in }\fQ_\tau\right\} \\
&=\left\{(s_{\tau'})\in\prod_{\tau'\in\fQ_\tau}G_{\tau'}:
\eta^{\tau'}_\tau(s_{\tau'})=s_\tau
\text{ for }\tau'\text{ in }\fQ_\tau\right\}
\end{align*}
which is locally compact by \cite[Proposition 2.1]{ilies2}.  
We let $\taunq$ denote the coarsest topology
which makes the map 
$s\mapsto\bigl(\eta_{\tau'}(s)\bigr):G\to G_{\fQ_\tau}$,
continuous.  We have that $\taunq\supseteq\tau'$ for every
$\tau'$ in $\fQ_\tau$.    

We now show that $\taunq\in\fQ_\tau$.  
First observe that $\eta^{\taunq}_\tau:G_{\taunq}\cong G_{\fQ_\tau}\to G_\tau$ 
is given by the map $(s_{\tau'})\mapsto s_\tau$.  With this identification
we have that
\[
\ker\eta^{\taunq}_\tau
=\underset{{\tau'\in\fQ_\tau}}{\underset{\longleftarrow}{\lim}}
\ker\eta^{\tau'}_\tau
\]
and is thus compact.  Moreover $\eta^{\taunq}_\tau$ is open, since
for any basic open set
\[
V=\prod_{\tau'\in\fQ_\tau\diagdown\{\tau'_1,\dots,\tau'_n\}}
G_{\tau'}\;\times\;\prod_{j=1}^n U_{\tau'_j}
\subset \prod_{\tau'\in\fQ_\tau}G_{\tau'}
\]
where each $U_{\tau_j'}$ is open in $G_{\tau_j'}$, we have
\[
\eta^{\taunq}_\tau(V\cap G_{\fQ_\tau})=
\bigcap_{j=1}^n\eta^{\tau_j'}_\tau(U_{\tau'_j})
\]
where each $\eta^{\tau_j'}_\tau(U_{\tau'_j})$ is open by assumption.

Finally, if there were $\tau_1$ in $\tg$ of which
$\taunq$ is a quotient, then by the first isomorphism theorem we would
have that $\tau_1\in\fQ_\tau$.  Hence $\tau_1\subseteq\taunq$.
Thus $\taunq$ is a non-quotient topology, and the unique such one
of which $\tau$ is a quotient. \endpf

For any locally compact group $G$ for which $\taunq=\tau\vee\tauap$
for any $\tau$ in $\tg$,
we obtained in \cite[Section 2.4]{ilies2} that {\it $\tnqg=\tg\vee\tauap$
and is thus an ideal in, and hence a subsemilattice of, 
the semilattice $(\tg,\vee)$.}  We note that for any $\tau$ for
which $G_\tau$ is maximally almost periodic, we have $\taunq=\tau\vee\tauap$.
Unfortunately, it is not clear whether
$\tnqg$ is a subsemilattice of $\tg$, in general.  However the
following is immediate.

\begin{corollary}\label{cor:tildevee}
If $\tau_1,\tau_2\in\tnqg$ define 
\[
\tau_1\tilde{\vee}\tau_2=(\tau_1\vee\tau_2)_{\mathrm{nq}}.
\]
Then $(\tnqg,\tilde{\vee})$ is a quotient semilattice of 
$(\tg,\vee)$.
\end{corollary}

The only inobvious aspect of this corollary is the associativity of $\tilde{\vee}$.
We observe that $\tau_1\vee(\tau_2\tilde{\vee}\tau_3)$ admits $\tau_1\vee\tau_2\vee\tau_3$
as a quotient, and hence, by Lemma \ref{ref:thelemma} (ii), is itself a quotient of 
$(\tau_1\vee\tau_2\vee\tau_3)_{nq}$.  Symmetrically, the same is true of 
$(\tau_1\tilde{\vee}\tau_2)\vee\tau_3$.  Hence
\[
\tau_1\tilde{\vee}(\tau_2\tilde{\vee}\tau_3)=(\tau_1\vee\tau_2\vee\tau_3)_{nq}
=(\tau_1\tilde{\vee}\tau_2)\tilde{\vee}\tau_3.
\]

Unless it can be shown that $\tau_1\tilde{\vee}\tau_2=\tau_1\vee\tau_2$
for all $\tau_1,\tau_2\in\tnqg$, some changes have to be made to the 
exposition in \cite[Section 4]{ilies2}, where $\vee$ must always be
replaced by, or understood to be, $\tilde{\vee}$.  Fortunately, this
change does not appear to affect any of the results or
proofs in this section, 
or any later part of the paper, in more than a cosmetic manner.

\section{Topology of the spine compactification}

As shown in \cite[Theorem 4.1]{ilies2}, with appropriate notational
changes as suggested by Corollary \ref{cor:tildevee},
the spectrum of the algebra $\spineg$ is given by
\[
G^*=\bigsqcup_{\fS\in\hdg}G_\fS
\]
where each $G_\fS$ is a projective limit over a hereditary
directed subset $\fS$ of $(\tnqg,\tilde{\vee})$.
The statement of \cite[Theorem 4.2]{ilies2} is flawed, and
should be replaced with the following.

\begin{theorem}\label{theo:topology}
The topology on $G^*$ is given as follows:  for any $s_0$ in $G^*$,
say $s_0\in G_{\fS_0}$ for some $\fS_0$ in $\hdg$, a neighbourhood
base at $s_0$ is formed by the sets
\begin{align*}
U(V_\tau;W_{\tau_1},\dots,W_{\tau_n})
=\bigm\{s\in G^*: & s\in G_\fS\text{ for some }\fS\supseteq\fS_\tau
\text{ in }\hdg  \\
&\text{with }s_\tau\in V_\tau,
\text{ and }s_{\tau_j}\in W_{\tau_j}\text{ if } \\
&\fS\supseteq\fS_{\tau_j},
\text{ for }j=1,\dots,n\bigm\}
\end{align*}
where $\tau\in\fS_0$, $\tau_1,\dots,\tau_n\in\tnqg\setminus\fS_0$, $V_\tau$
is an open neighbourhood of $s_{0,\tau}$ in $G_\tau$, and each
$W_{\tau_j}$ is a cocompact subset of $G_{\tau_j}$.
\end{theorem}

Pekka Salmi has pointed out to us that the error in the description of 
\cite[Theorem 4.2]{ilies2}, implies that all groups $G_\fS$, for $\fS\in\hdg$
are locally compact.  This is false, as is implicit in 
\cite[Section 6.3]{ilies2}, or is shown in \cite[Theorem 2]{dunklr}.

Unfortunately, the proof of \cite[Theorem 4.2]{ilies2} requires slight 
modification.

\medskip
\noindent {\bf Proof of Theorem \ref{theo:topology}.}
We should first observe that the family of sets described above
indeed is a base for a topology.  It is straightforward to check
that for $\tau,\tau'\in\fS_0$
and $\tau_1,\dots,\tau_n,\tau_1',\dots,\tau_m'$ 
in $\tnqg\setminus\fS_0$ that
\begin{align*}
U&(V_\tau;W_{\tau_1},\dots,W_{\tau_n})\cap
U(V_{\tau'};W_{\tau_1'},\dots,W_{\tau_m'}) \\
&\supset U\left((\eta^{\tau\tilde{\vee}\tau'}_\tau)^{-1}(V_\tau)\cap
(\eta^{\tau\tilde{\vee}\tau'}_{\tau'})^{-1}(V_{\tau'});
W_{\tau_1},\dots,W_{\tau_n},W_{\tau_1'},\dots,W_{\tau_n'}\right)
\end{align*}
for neighbourhoods $V_\tau$ of $s_{0,\tau}$, $V_{\tau'}$ of $s_{0,\tau'}$,
and cocompact subsets $W_{\tau_1},\dots,W_{\tau_m'}$ of $G_{\tau_1},\dots,
G_{\tau_m'}$, respectively.

We note that for a net $(s_i)_{i\in I}$ in $G^*$, the following are equivalent:

(i) $s_i\to s_0$ in $G^*$ with the topology described above;

(ii) for each $\tau_0\in\fS_0$ there is $i_0$ such that $i\geq i_0$
implies that $s_i\in G_{\fS_i}$ for some $\fS_i\supset\fS_{\tau_0}$ and
$\lim_{i\geq i_0}s_{i,\tau}=s_{0,\tau}$; and for each $\tau\not\in\fS_0$
for which $I_\tau=\{i:s_i\in\fS_i\text{ for some }\fS_i\supset\fS_\tau\}$
is admits no maximal element in $I$, and for any co-compact $W_\tau\subset G_\tau$,
there is $i_\tau$ in $I$ 
for which $s_{i,\tau}\in W_\tau$ if $i\in I_\tau$ and $i\geq i_\tau$;

(iii) $\chi_{s_i}\to \chi_{s_0}$ weak* in $\spineg^*$.

The equivalence of (i) and (ii) is clear.  If we write $u$ in $\spineg$
as $u=\sum_{\tau\in\tnqg}u_\tau$ as in \cite[(4.2)]{ilies2}, and suppose
(ii), above, then
\[
\chi_{s_i}(u)=\sum_{\tau\in\fS_0\setminus\fS_i}\hat{u}_\tau(s_{i,\tau})+
\sum_{\tau\in\fS_i\setminus\fS_0}\hat{u}_\tau(s_{i,\tau})\to 
\sum_{\tau\in\fS_0}\hat{u}_\tau(s_{0,\tau})=\chi_{s_0}(u)
\]
where $s_i\in\fS_i$ for each $i$; which shows (iii).  Likewise,
selecting $u=u_\tau$, a repeat of the computation above shows that
(iii) implies (ii). \endpf

We remark that
\cite[Corollary 4.3]{ilies2} remains unchanged.  The neighbourhood in the proof
of part (ii), therein, should be changed to
\[
U(V_{\tau_0})=\{s\in G^*:s\in G_{\fS'}\text{ for some }\fS'\supset\fS_{\tau_0}
\text{ and }S_{\tau_0}\in V_{\tau_0}\}.
\]
Furthermore, a modest change must be made in the proof of 
\cite[Proposition 4.6 (ii)]{ilies2}.  
If a net $(e_i)$ of idempotents converges to
$s$ in $G_\fS$, then for all $\tau\in\fS$, there is an $i_\tau$ for which
$\fS_i\supset\fS_\tau$, where $e_i\in G_{\fS_i}$, for $i\geq i_\tau$, and
$\lim_{i\geq i_\tau}e_{i,\tau}=s_\tau$; and for $\tau$ in $\tnqg\setminus\fS$,
there is $i_\tau$ for which $\fS_i\not\supset\fS_\tau$ for $i\geq i_\tau$,
which may be seen by observing that otherwise such $e_{i,\tau}$ must be within
the cocompact set $W_\tau=G_\tau\setminus\{\eta_\tau(e)\}$.

\section{On abelian groups}

The fourth paragraph of the proof of \cite[Theorem 5.1]{ilies2}
contains an error in its claim that the map
from $\widehat{G_\tau}$ to the diagonal subgroup of 
$\widehat{G_\tau}\times \widehat{G}_d$ is bicontinuous is false, for 
obviously only its inverse is continuous.  Fortunately the claim which 
that paragraph is attempting to establish, namely that
{\em if $\tau$ in $\tg$ is such that $\widehat{G_\tau}$ is the group 
$\widehat{G}$ but with a finer locally compact group topology $\hat{\tau}$, 
then $\tau\in\tnqg$}, remains true.  Recall that for any $\tau'$ in $\tg$, the dual 
map $\widehat{\eta_{\tau'}}:\widehat{G_{\tau'}}\to\widehat{G}$
is continuous and injective.  The assumptions on $\tau$ above imply  that
$\widehat{\eta_\tau}$ is also surjective.  
Since the map $\eta^{\taunq}_\tau:G_{\taunq}\to G_\tau$ is proper, its
dual map $\widehat{\eta^{\taunq}_\tau}$ 
is open; see \cite[(23.24)(d)]{hewittrI}, for example.  
Since $\eta_\tau=\eta^{\taunq}_{\tau}\comp\eta_{\taunq}$ we have
$\widehat{\eta_\tau}=\widehat{\eta_{\taunq}}\comp\widehat{\eta^{\taunq}_{\tau}}$.
Thus it follows that $\widehat{\eta^{\taunq}_{\tau}}$ is continuous, bijective
and open, hence $\hat{\tau}=\widehat{\taunq}$.  Thus by Pontryagin duality
$\tau=\taunq$.



\medskip
{\sc  Monica Ilie, Department of Mathematical Sciences, Lakehead University,
955 Oliver Road,  Thunder Bay, ON, P7B\,5E1, Canada} \hfill
{\tt milie@lakeheadu.ca}

\medskip
{\sc Department of Pure Mathematics, University of Waterloo,
Waterloo, ON, N2L 3G1, Canada} \hfill {\tt nspronk@uwaterloo.ca}

\end{document}